\documentclass{amsart}
\usepackage{amssymb}
\usepackage{amsthm}
\usepackage{amstext}
\usepackage{amsbsy}
\usepackage{amscd}
\usepackage{enumerate}

\newtheorem{thm}{Theorem}[section]

\newtheorem{quest}[thm]{Question}

\newtheorem{defn}[thm]{Definition}

\newtheorem*{thma}{Theorem A}
\newtheorem*{thmb}{Theorem B}

\newtheorem*{thmh}{Theorem H}
\newtheorem*{thmg}{Theorem G}

\newcommand{\R}{\mathbb{R}}
\newcommand{\Q}{\mathbb{Q}}

\newcommand{\Z}{\mathbb{Z}}

\newcommand{\dist}{\mathop{\rm{dist}}}
\newcommand{\cS}{{\mathcal{S}}}
\def\zero{{(0)}}

\allowdisplaybreaks

\begin{document}

\title[The spectrum of a q.p.\ Schr\"odinger operator is homogeneous]{The spectrum of a Schr\"odinger operator with small quasi-periodic potential is homogeneous}

\author[D. Damanik]{David Damanik}
\address{Department of Mathematics, Rice University, Houston TX 77005, U.S.A.}
\email{damanik@rice.edu}

\author[M. Lukic]{Milivoje Lukic}
\address{Department of Mathematics, Rice University, Houston TX 77005, U.S.A.}
\email{milivoje.lukic@rice.edu}

\author[M. Goldstein]{Michael Goldstein}
\address{Department of Mathematics, University of Toronto, Bahen Centre, 40 St. George St., Toronto, Ontario, CANADA M5S 2E4}
\email{gold@math.toronto.edu}

\thanks{D.~D.\ was partially supported by NSF grants DMS--1067988 and DMS--1361625. M.~G.\ was partially supported by NSERC. M.~G.\ expresses his gratitude for the hospitality during a stay at the Institute of Mathematics at the University of Stony Brook in May 2014. M.~L.\ was partially supported by NSF grant DMS--1301582.}

\date{\today}

\begin{abstract}
We consider the quasi-periodic Schr\"odinger operator
$$
[H \psi](x) = -\psi''(x) + V(x) \psi(x)
$$
in $L^2(\R)$, where the potential is given by
$$
V(x) = \sum_{m \in \Z^\nu \setminus \{ 0 \}} c(m)\exp (2\pi i m \omega x)
$$
with a Diophantine frequency vector $\omega = (\omega_1, \dots, \omega_\nu) \in \R^\nu$ and exponentially decaying Fourier coefficients $|c(m)| \le \varepsilon \exp(-\kappa_0|m|)$. In the regime of small $\varepsilon > 0$ we show that the spectrum of the operator $H$ is homogeneous in the sense of Carleson.
\end{abstract}

\maketitle

\sloppy

\section{Introduction and Statement of the Main Result}\label{sec.1}

Consider a Schr\"odinger operator
\begin{equation}\label{eq:so}
[H \psi](x) = - \psi''(x) +  V(x) \psi(x)
\end{equation}
in $L^2(\R)$. Associated with such an operator is spectral information, such as the spectrum and the spectral measures.

For the most part, the spectral analysis of operators of this kind breaks into two branches, namely direct spectral analysis, where the potential is given and one seeks information about the spectrum and/or the spectral measures, and inverse spectral theory, where information about the spectrum and/or the spectral measures is given and one seeks information about the potential.

In very special cases one can obtain two-way results of this kind, where certain classes of potentials are in one-to-one correspondence with certain classes of spectra and/or spectral measures. Results of this special nature are called ``gems of spectral theory'' in Barry Simon's monograph \cite{Simon}.

Most results in the spectral theory of Schr\"odinger (or related) operators, however, are not of this special nature, and they are strictly one-way results. That is, for a certain class of potentials one can prove certain spectral properties, or conversely, for a certain class of potentials defined via the spectral properties of the associated operators, one can prove certain statements. Whenever this situation arises, there is a gap in our understanding of the spectral problem at hand, and we do not have a complete characterization of the class of potentials or spectral features that correspond to the class on the other side under consideration. In this case there is a natural interest in closing this gap.

\medskip

After these general introductory remarks, let us be more concrete. Assuming that $V$ is almost periodic (i.e., the set of its translates is relatively compact in the uniform topology), there is extensive literature on the direct spectral problem, that is, proving statements about the spectrum of the operator $H$ and the type of the spectral measures associated with it. We direct the reader to the recent surveys \cite{D, JM} and the references therein. The appearance of Cantor sets as spectra turns out to be typical, and the spectral measures can be of all possible types, but essentially they have a tendency to be purely absolutely continuous for small potentials or for large energies, and pure point for large potentials and small energies, provided that the potential has sufficient regularity properties (in the highly irregular case, the appearance of purely singular continuous spectral measures is typical). The case of regular small quasi-periodic potentials is quite well understood; see, for example, \cite{DG, E}. Here, indeed the spectrum is a Cantor set and the spectral measures are purely absolutely continuous. It is known due to work of Kotani \cite{K84, K97} and Remling \cite{R07} that, as a consequence of absolute continuity, the operators in question are reflectionless, that is, the boundary values of the diagonal elements of the Green function are purely imaginary Lebesgue almost everywhere on the (absolutely continuous) spectrum.

In the converse direction, there has been extensive work on reflectionless Schr\"odinger operators for certain prescribed spectra. Here one fixes a set $\cS$ and considers the set of all potentials $V$ such that the associated Schr\"odinger operator has spectrum $\cS$ and is reflectionless on it. The goal is then to find out as much as possible about this class of potentials. Fundamental work in this direction can be found, for example, in  \cite{Cr, SY}. It is natural to ask about conditions on $\cS$ that ensure that all the potentials associated with it are almost periodic and all spectral measures are purely absolutely continuous. It turns out that the following condition (cf.~\cite{Ca}) does the job.

\begin{defn}
A closed set
$$
\mathcal{E} = \mathbb{R} \setminus \bigcup_n (E^-_n,E^+_n)
$$
is called homogeneous if there is $\tau > 0$ such that for any $E \in \mathcal{E}$ and any $\sigma > 0$, we have $|(E - \sigma, E + \sigma) \cap \mathcal{E}| > \tau \sigma$.
\end{defn}

Assuming finite total gap length,\footnote{Finite total gap length means that the sum of the lengths of the bounded gaps of the spectrum is finite.} it was shown by Sodin and Yuditskii \cite{SY} that the homogeneity of $\cS$ implies the almost periodicity of the associated potentials, and it was shown by Gesztesy and Yuditskii \cite{GY} that the homogeneity of $\cS$ implies the absolute continuity of the associated spectral measures (see also the more recent paper \cite{PR}). On the other hand, it is known that neither consequence holds without a suitable assumption on $\cS$, such as for example homogeneity. Namely, Poltoratski and Remling study conditions on the set related to the presence of associated reflectionless measures with non-trivial singular component \cite{PR} and, working out the continuum analog of work by Volberg and Yuditskii \cite{VY}, Damanik and Yuditskii showed that there are sets $\cS$ such that all associated potentials are not almost periodic \cite{DY}.\footnote{The sets in \cite{DY} are essentially explicit. With a non-explicit set, the statement can also be derived from \cite{A4}.}

Apparently, one has a good way of passing to almost periodicity and absolute continuity from either side, that is, one has conditions on potentials that ensure these properties, and one has spectral conditions that ensure these properties as well. Can one find a link between these two rather different sets of results? It is the main purpose of this paper to establish such a link. Namely, we will show that the operators studied in \cite{DG}, which is a direct spectral analysis, have homogeneous spectra. Since they also clearly have finite total gap length (and the operators are reflectionless as pointed out above), this puts the operators studied in \cite{DG} inside the scope of the relevant literature on inverse spectral analysis, and in particular the papers \cite{GY, SY}. We mention in passing that another link of this nature will be established in \cite{DGL1, DGL2}. Namely it is shown there that if $\cS$ is of the type that arises in \cite{DG}, then all potentials for which the associated operator is reflectionless are again of the type studied in \cite{DG}. In particular, they are all regular and quasi-periodic with the same frequency vector. That is, for these sets $\cS$, the potentials studied in the inverse spectral theory approach are put inside the scope of the work \cite{DG} on the direct spectral problem.

\bigskip

Let us now proceed to the statement of the main result of this paper. Let $U(\theta)$ be a real function on the torus $\mathbb{T}^\nu$,
$$
U(\theta) = \sum_{n \in \mathbb{Z}^\nu \setminus \{ 0 \}} c(n) e^{2 \pi i n\theta}\ , \quad \theta \in \mathbb{T}^\nu.
$$

Let $\omega = (\omega_1, \dots, \omega_\nu) \in \mathbb{R}^\nu$. Assume that the following Diophantine condition holds,
\begin{equation}\label{eq:1PAI7-5-85a}
|n \omega| \ge a_0 |n|^{-b_0}, \quad n \in \mathbb{Z}^\nu \setminus \{ 0 \}
\end{equation}
for some
$$
0 < a_0 < 1,\quad \nu < b_0 < \infty.
$$

Let $V(x) = U(x \omega)$ and consider the Schr\"odinger operator \eqref{eq:so}. Assume that $U$ is real-analytic, that is, the Fourier coefficients $c(n)$ obey
\begin{align*}
\overline{c(n)} & = c(-n), \quad n \in \Z^\nu \setminus \{ 0 \}, \\
|c(n)| & \le  \varepsilon\exp(-\kappa_0|n|), \quad n \in \Z^\nu \setminus \{ 0 \},
\end{align*}
with $\varepsilon > 0$, $0 < \kappa_0 \le 1$.

Our main result reads as follows:

\begin{thmh}
There exists $\varepsilon_0 = \varepsilon_0(\kappa_0, a_0, b_0) > 0$ such that for $0 < \varepsilon < \varepsilon_0$, the spectrum of the operator $H$ is homogeneous with $\tau = 1/2$.
\end{thmh}

The homogeneity of the spectrum is a consequence of detailed quantitative results we can establish for the structure of the gaps of the spectrum. Since the latter results are of independent interest, we state them separately in the following theorem.

\begin{thmg}
There exists $\varepsilon_0 = \varepsilon_0(\kappa_0, a_0, b_0) > 0$ such that for $0 < \varepsilon < \varepsilon_0$, the gaps in spectrum of the operator $H$ can be labeled as $G_m = (E_m^-, E_m^+)$, $m \in \mathbb{Z}^\nu \setminus \{ 0 \}$, $G_0 = (-\infty, \underline{E})$ so that the following conditions hold:
\begin{enumerate}[{\rm (i)}]

\item For every $m \in \mathbb{Z}^\nu \setminus \{ 0 \}$, we have
$$
E^+_m - E^-_m \le 2 \varepsilon \exp \Big( -\frac{\kappa_0}{2} |m| \Big).
$$

\item For every $m, m' \in \mathbb{Z}^\nu \setminus \{ 0 \}$ with $m' \neq m$ and $|m'| \ge |m|$, we have
$$
\dist ([E_m^-, E_m^+], [E_{m'}^-, E_{m'}^+])) \ge a |m'|^{-b},
$$
where $a, b > 0$ are constants depending on $a_0, b_0, \kappa_0, \nu$.

\item For every $m \in \mathbb{Z}^\nu \setminus \{ 0 \}$, we have
$$
E_m^- - \underline{E}\ge a|m|^{-b}.
$$

\end{enumerate}
\end{thmg}

In the setting described above, Damanik and Goldstein established in \cite{DG} a rather detailed description of the spectrum and the generalized eigenfunctions, which turn out to be of Floquet type. As a consequence of this description and the work of Kotani and Remling mentioned above \cite{K84, K97, R07} it follows that the operator $H$ is indeed reflectionless on its spectrum, which in turn has finite total gap length by the estimates in \cite{DG}. Thus, the only missing piece for one to apply the theory of Gesztesy, Sodin and Yuditskii \cite{GY, SY} is the homogeneity of the spectrum. This missing piece is established by Theorem~H. On the more conceptual level discussed earlier, Theorem~H therefore provides a link between the direct and the inverse spectral theory approach to almost periodicity and absolute continuity.

\bigskip

The structure of the remainder of the paper is as follows. We prove Theorems~H and G in Section~\ref{sec.2} and discuss natural questions for further study in Section~\ref{sec.3}.

\section{Proof of Theorems H and G}\label{sec.2}

In this section we prove Theorem~H and G. The main results from \cite{DG} play a crucial role in these proofs. These results, given in Theorems~A and B in \cite{DG} and restated below, describe the spectrum and the generalized eigenfunctions of the operator $H$ with a small analytic quasi-periodic potential $V$ and establish a two-way connection between the decay of the Fourier coefficients of $V$ and the size of the gaps of the spectrum of $H$.

Let us recall these results from \cite{DG}. Set
\begin{align*}
k_n & = -n\omega/2, \quad n \in \Z^\nu \setminus \{0\}, \quad \mathcal{K}(\omega) = \{ k_n : n \in \Z^\nu \setminus \{0\} \}, \\
\mathfrak{J}_n & = ( k_n - \delta(n), k_n + \delta(n) ), \quad \delta(n) = a_0 (1 + |n|)^{-b_0-3}, \quad n \in \Z^\nu \setminus \{0\}, \\
\mathfrak{R}(k) & = \{ n \in \Z^\nu \setminus \{0\} : k \in \mathfrak{J}_n \}, \quad \mathfrak{G} = \{ k : |\mathfrak{R}(k)| < \infty \},
\end{align*}
where $a_0,b_0$ are as in the Diophantine condition \eqref{eq:1PAI7-5-85a}. Let $k \in \mathfrak{G}$ be such that $|\mathfrak{R}(k)| > 0$. Due to the Diophantine condition, one can enumerate the points of $\mathfrak{R}(k)$ as $n^{(\ell)}(k)$, $\ell = 0, \dots, \ell(k)$, $1 + \ell(k) = |\mathfrak{R}(k)|$, so that $|n^{(\ell)}(k)| < |n^{(\ell+1)}(k)|$. Set
\begin{align*}
T_{m}(n) & = m - n ,\quad m, n \in \mathbb{Z}^\nu, \\
\mathfrak{m}^{(0)}(k) & = \{ 0, n^{(0)}(k) \}, \\
\mathfrak{m}^{(\ell)}(k) & = \mathfrak{m}^{(\ell-1)}(k) \cup T_{n^{(\ell)}(k)}(\mathfrak{m}^{(\ell-1)}(k)), \quad \ell = 1, \dots, \ell(k).
\end{align*}

The following pair of theorems was established in \cite{DG}.

\begin{thma}
There exists $\varepsilon_0 = \varepsilon_0(\kappa_0, a_0, b_0) > 0$ such that for $0 < \varepsilon < \varepsilon_0$ and $k \in \mathfrak{G} \setminus \frac{\omega}{2}(\Z^\nu \setminus \{0\})$, there exist $E(k) \in \mathbb{R}$ and $\varphi(k) := (\varphi(n;k))_{n \in \Z^\nu}$ such that the following conditions hold:

$(a)$ $\varphi(0; k) = 1$,
\begin{align*}
|\varphi(n ;k)| & \le \varepsilon^{1/2} \sum_{m \in \mathfrak{m}^{(\ell)}} \exp \Big( -\frac{7}{8} \kappa_0 |n-m| \Big), \quad \text{ $n \notin \mathfrak{m}^{(\ell(k))}(k)$}, \\
|\varphi(m;k)| & \le 2 \quad \text{for any $m \in \mathfrak{m}^{(\ell(k))}(k)$.}
\end{align*}

$(b)$ The function
$$
\psi(k, x) = \sum_{n \in \Z^\nu} \varphi(n; k) e^{2 \pi i x (n \omega + k)}
$$
is well-defined and obeys
$$
- \psi''(k,x) + V(x) \psi(k,x) = E(k) \psi(k,x).
$$

$(c)$
$$
E(k) = E(-k), \quad \varphi(n;-k) = \overline{\varphi(-n; k)}, \quad \psi(-k, x)=\overline{\psi(k, x)},
$$
\begin{equation}\label{eq:6Ekk1EGT11}
\begin{split}
(k^0)^2 (k - k_1)^2  < E(k) - E(k_1) < 2k (k - k_1) + 2 \varepsilon \sum_{k_1 < k_{n} < k} \delta(n), \quad 0 < k - k_1 < 1/4, \; k_1 > 0,
\end{split}
\end{equation}
where $k^\zero := \min(\varepsilon_0, k/1024)$.

$(d)$ The spectrum of $H$ consists of the following set,
$$
\cS = [E(0) , \infty) \setminus \bigcup_{m \in \Z^\nu \setminus \{0\} : E^-( k_m) < E^+( k_m)} (E^-( k_m), E^+( k_m)),
$$
where
$$
E^\pm(k_m) = \lim_{k \to k_m \pm 0, \; k \in \mathfrak{G} \setminus \mathcal{K}(\omega)} E(k), \quad \text{ for $k_m>0$.}
$$
\end{thma}

\begin{thmb}
{\rm (a)} The gaps $(E^-(k_m), E^+(k_m))$ in Theorem~A obey $E^+(k_m) - E^-(k_m) \le 2 \varepsilon \exp(-\frac{\kappa_0}{2} |m|)$.

 {\rm (b)} Using the notation from Theorem~A, there exists $\varepsilon^\zero > 0$ such that if the gaps $(E^-(k_m), E^+(k_m))$ obey $E^+(k_m) - E^-(k_m) \le \varepsilon \exp(-\kappa |m|)$ with $0 < \varepsilon < \varepsilon^\zero$, $\kappa > 4 \kappa_0$, then, in fact, the Fourier coefficients $c(m)$ obey $|c(m)| \le \varepsilon^{1/2} \exp(-\frac{\kappa}{2} |m|)$.
\end{thmb}

We are now in position to prove Theorems~G and H.

\begin{proof}[Proof of Theorem~$G$]
Consider the $\varepsilon_0 = \varepsilon_0(\kappa_0, a_0, b_0) > 0$ from Theorem~A and label the gaps as in part (d) of Theorem~A.

(i) This statement follows from part (a) of Theorem~B.

(ii) Recall  that
$$
|m\omega| > a_0 |m|^{-b_0}, m\neq 0.
$$
In what follows we denote by $a_j$ constants depending on $a_0, b_0, \kappa_0, \nu$. Let $m' \neq m$, $|m'| \ge |m|$ be arbitrary. Then,
$$
|k_m - k_{m'}| = |(m-m')\omega|/2 \ge a_0 (2 |m'|)^{-b_0}/2 \ge a_1 |m'|^{-b_0}.
$$
Assume for instance that $k_{m'} > k_{m} > 0$. Due to \eqref{eq:6Ekk1EGT11} in Theorem~A, we have
$$
E^-(k_{m'}) - E^+(k_{m}) > (k^0)^2 (k_{m'} - k_m)^2 \ge (k^0)^2 a_2 |m'|^{-2b_0},
$$
where
$$
k^\zero := \min (\varepsilon_0, k_{m} /1024) \ge \varepsilon_0 a_3 |m'|^{-b_0}.
$$
Thus,
$$
E^-(k_{m'}) - E^+(k_{m}) > (k^0)^2 (k_{m'} - k_m)^2 \ge \varepsilon_0^2 a_4 |m'|^{-4 b_0} = a_5 |m'|^{-4 b_0},
$$
which means that
$$
\dist ([E_m^-, E_m^+], [E_{m'}^-, E_{m'}^+]) \ge a |m'|^{-b}.
$$
The remaining cases are completely similar.

(iii) The proof of this statement is completely similar to the proof of (ii) and we omit it.
\end{proof}

\begin{proof}[Proof of Theorem~$H$]
Let $E \in \mathcal{S}$, $\sigma > 0$. Set
$$
\mathfrak{C}(E,\sigma)=\{m\neq 0: (E^-_m,E^+_m)\cap (E-\sigma,E+\sigma)\neq \emptyset\}.
$$
Using the notation from Theorem $A$, assume first that
$$
(-\infty, \underline{E}) \cap (E - \sigma, E + \sigma) = \emptyset.
$$
Pick $m_0 = m_0(E, \sigma)$ so that $|m_0| = \min_{m \in \mathfrak{C}(E, \sigma)}|m|$. Note that for any $m \in \mathfrak{C}(E,\sigma)$, we have
$$
\dist ([E_m^-, E_m^+], [E_{m_0}^-, E_{m_0}^+]) \le 2 \sigma.
$$
On the other hand, by part (ii) of Theorem~G,
$$
\dist ([E_m^-, E_m^+], [E_{m_0}^-, E_{m_0}^+])\ge a|m|^{-b},
$$
where $a, b > 0$ are constants depending on $a_0, b_0, \kappa_0, \nu$. Therefore,
$$
|m| \ge \alpha \sigma^{-\beta},
$$
where $\alpha, \beta > 0$ are constants depending on $a_0, b_0, \kappa_0, \nu$. Due to part (i) of Theorem~G,
\begin{equation}\label{eq:3gapslenthsrep1}
E^+_{m} - E^-_{m} < 2 \varepsilon \exp \Big( -\frac{\kappa_0}{2} |m| \Big).
\end{equation}
Thus,
\begin{align*}
\sum_{m \in \mathfrak{C}(E, \sigma) \setminus \{ m_0 \}} \big| (E^-_{m}, E^+_{m}) \cap (E - \sigma, E + \sigma) \big| & \le \sum_{m \in \mathfrak{C}(E, \sigma) \setminus \{ m_0 \}} E^+_m - E^-_m \\
& \le 2 \varepsilon \sum_{|m| \ge \alpha \sigma^{-\beta}} \exp \Big( -\frac{\kappa_0}{2} |m| \Big) \\
& < \sigma/2,
\end{align*}
provided $\sigma \le \sigma_0(a_0, b_0, \kappa_0, \nu)$. Note that since $E \in \mathcal{S}$, $E \notin (E_{m_0}^-, E_{m_0}^+)$. Hence,
$$
\big| (E^-_{m_0}, E^+_{m_0}) \cap (E - \sigma, E + \sigma) \big| \le \sigma.
$$
Thus,
\begin{align*}
\big| (E - \sigma, E + \sigma) \cap \mathcal{S} \big| & \ge 2 \sigma - \big| (-\infty, \underline{E}) \cap (E - \sigma, E + \sigma) \big| - \big| (E^-_{m_0}, E^+_{m_0}) \cap (E - \sigma, E + \sigma) \big| \\
& \qquad - \sum_{m \in \mathfrak{C}(E, \sigma) \setminus \{ m_0 \}} \big| (E^-_{m}, E^+_{m}) \cap (E - \sigma, E + \sigma) \big| \\
& > 2 \sigma - 0 - \sigma - \sigma/2 \\
& = \sigma/2,
\end{align*}
provided $\sigma \le \sigma_0(a_0, b_0, \kappa_0, \nu)$.

Now assume that
$$
(-\infty, \underline{E}) \cap (E - \sigma, E + \sigma) \neq \emptyset.
$$
Note that for any $m \in \mathfrak{C}(E, \sigma)$, we have
$$
E_m^- - \underline{E} \le 2 \sigma.
$$
On the other hand, by part (iii) of Theorem~G, we have
$$
E_m^- - \underline{E} \ge a |m|^{-b},
$$
where $a, b > 0$ are constants depending on $a_0, b_0, \kappa_0, \nu$. Therefore,
$$
|m| \ge \alpha \sigma^{-\beta},
$$
where $\alpha, \beta > 0$ are constants depending on $a_0, b_0, \kappa_0, \nu$. Just as above we may conclude that
$$
\sum_{m \in \mathfrak{C}(E, \sigma)} \big| (E^-_{m}, E^+_{m}) \cap (E - \sigma, E + \sigma) \big| < \sigma/2,
$$
provided $\sigma \le \sigma_0(a_0, b_0, \kappa_0, \nu)$.  Note that since $E \in \mathcal{S}$, $E \notin (-\infty, \underline{E})$. Hence,
$$
\big| (-\infty, \underline{E}) \cap (E - \sigma, E + \sigma) \big| \le \sigma.
$$

Due to \eqref{eq:3gapslenthsrep1},
$$
\sum_m (E^+_{m} - E^-_{m}) < C(\kappa_0, \nu) \varepsilon < \sigma_0(a_0, b_0, \kappa_0, \nu)/2,
$$
provided $0 < \varepsilon < \varepsilon_0 = \varepsilon_0(\kappa_0, a_0, b_0) > 0$. Therefore, for any interval $(E - \sigma, E + \sigma)$ with $E \in \mathcal{S}$ and $\sigma > \sigma_0(a_0, b_0, \kappa_0, \nu)$, we have
\begin{align*}
\big| (E - \sigma, E + \sigma) \cap \mathcal{S} \big| & \ge 2 \sigma - \big| (-\infty, \underline{E}) \cap (E - \sigma, E + \sigma) \big| - \sum_{m} \big| (E^-_{m}, E^+_{m}) \big| \\
& > 2 \sigma - \sigma - \sigma/2 \\
& = \sigma/2,
\end{align*}
which shows that the desired estimate holds in the second case as well. This concludes the proof.
\end{proof}

\section{Some Remarks and Open Problems}\label{sec.3}

We have seen that Schr\"odinger operators with small analytic quasi-periodic potentials have homogeneous spectrum. It is natural to ask whether this property persists if $\varepsilon$ is increased somewhat. It is known that, as $\varepsilon$ is increased, the pure absolute continuity of the spectrum does not persist. Indeed for large enough coupling, the spectrum will be pure point with exponentially localized eigenfunctions (i.e., Anderson localization holds) in the lower energy region. This is particularly well understood for the discrete counterpart of this problem, but we expect very strongly that the continuum versions of the statements known in the discrete case indeed do hold. Moreover, judging again by the analogy with the discrete case, the transition from absolute continuity to localization may well go through a critical regime at which the homogeneity of the spectrum breaks down.\footnote{Concretely, for the almost Mathieu operator this transition occurs at coupling $\lambda = 1$, and at this value of the coupling constant, the spectrum has zero Lebesgue measure.} It is of course of interest to compare the values of the coupling constant where we experience a breakdown of absolute continuity and homogeneity, respectively. Let us state one question in this spirit explicitly.

\begin{quest}\label{q.1}
Suppose that $\varepsilon_1 > 0$ is such that for $0 < \varepsilon < \varepsilon_1$, the spectrum of $H$ is purely absolutely continuous with generalized eigenfunctions of Floquet type for almost every energy in the spectrum. Is it true that for each $0 < \varepsilon < \varepsilon_1$, the spectrum of $H$ is homogeneous?
\end{quest}

Another interesting research direction is to explore the same issues in the discrete setting. Given that inverse spectral theory plays a role in this study, one needs to study the class of Jacobi matrices, that is, operators
$$
[J \psi]_n = a_n \psi_{n+1} + b_n \psi_n + a_{n-1} \psi_{n-1}
$$
in $\ell^2(\Z)$ with $a_n > 0$ and $b_n \in \R$. The inverse spectral theory aspects were worked out by Sodin and Yuditskii as well \cite{SY2}. The connection between absolute continuity and reflectionlessness follows from the work of Kotani \cite{K97} and Remling \cite{R11}. There is a very large number of results on the direct spectral problem for Jacobi matrices with almost periodic coefficients, we refer the reader again to the recent surveys \cite{D, JM} and the references therein. In fact, there are more results in the discrete setting than in the continuum setting. For example, Avila's global theory for analytic quasi-periodic one-frequency potentials \cite{A1, A2, A3} currently exists only in the discrete setting. This suggests that one should try to establish discrete versions of \cite{DG, DGL1, DGL2} and of the present paper.

By far the most heavily studied discrete Schr\"odinger operator with a quasi-periodic potential is the almost Mathieu operator, that is, the Jacobi matrix with $a_n = 1$ and $b_n = 2 \lambda \cos (2 \pi [n\omega + \theta])$ with $\lambda > 0$, $\omega \in \R \setminus \Q$ and $\theta \in \R$. It is known that the spectrum of this operator does not depend on $\theta$, and may therefore be denoted by $\sigma_{\lambda,\omega}$, and has Lebesgue measure $\mathrm{Leb}(\sigma_{\lambda,\omega}) = 4 | 1 - \lambda |$. Moreover, all spectral measures are purely absolutely continuous for $\lambda < 1$ and purely singular for $\lambda \ge 1$. There is a rather large number of papers contributing to these statements, compare \cite{D, JM}. In particular, the operator (has spectrum of positive Lebesgue measure and) is reflectionless when $\lambda < 1$. On the other hand it is not known whether there are any parameter values for which the spectrum of this operator is homogeneous. Given that almost everything about this operator is known due to extensive investigation over many decades, this absence of understanding is somewhat unsettling. Let us state this explicitly as a question.

\begin{quest}\label{q.2}
For which values of $\lambda$ is the spectrum of the almost Mathieu operator homogeneous?
\end{quest}

Working out the discrete analog of \cite{DG} and the present paper as suggested above could potentially solve this problem in the regime of small $\lambda$. We intend to address this in a forthcoming paper.

Note that by Aubry duality, the spectrum at $\lambda$ is homogeneous if and only if it is homogeneous at $1/\lambda$. Thus, we expect homogeneity to hold at least for $\lambda$ sufficiently small and sufficiently large.

However, homogeneity will not hold for all $\lambda > 0$. At the critical value $\lambda = 1$, the spectrum has zero Lebesgue measure and hence homogeneity fails for trivial reasons. It is far from clear at this point whether homogeneity will hold for all $\lambda < 1$ or whether it breaks down at some smaller threshold. In this context the work of Helffer and Sj\"ostrand \cite{HS} may be relevant. In that paper, the gap structure for the case of critical coupling $\lambda = 1$ and suitable frequencies $\omega$ was studied in great detail. If it is possible to extend their analysis to non-critical $\lambda$, this may allow one to address the homogeneity issue in the almost critical regime. This appears to be quite difficult, however.

\section*{Acknowledgment}

We thank Bernard Helffer, Svetlana Jitomirskaya and Johannes Sj\"ostrand for useful conversations about the homogeneity question for the almost Mathieu operator.


\begin{thebibliography}{DGL2}

\bibitem[A1]{A1} Avila, A. {\em Global theory of one-frequency Schr\"odinger operators I: stratifed analyticity of the Lyapunov exponent and the boundary of nonuniform hyperbolicity.} Preprint (arXiv:0905.3902).

\bibitem[A2]{A2} Avila, A. {\em Global theory of one-frequency Schr\"odinger operators II: acriticality and finiteness of phase transitions for typical potentials.} Preprint.

\bibitem[A3]{A3} Avila, A. {\em Almost reducibility and absolute continuity I.} Preprint (arXiv:1006.0704).

\bibitem[A4]{A4} Avila, A. {\em On the Kotani-Last and Schr\"odinger conjectures.} To appear in J.\ Amer.\ Math.\ Soc.\ (arXiv:1210.6325).

\bibitem[Ca]{Ca} Carleson, L. {\em On $H^\infty$ in multiply connected domains.} Harmonic Analysis. Conference in honor of Antony Zygmund.
 vol. II, (1983), 349--382.

\bibitem[Cr]{Cr} Craig, W. {\em Trace formula for the Schr\"{o}dinger operator on the line.} Commun.\ Math.\ Phys.\ \textbf{126} (1989), 379--407.

\bibitem[D]{D} Damanik, D. {\em Schr\"odinger operators with dynamically defined potentials: a survey.} In preparation.

\bibitem[DG]{DG} Damanik, D., Goldstein, M. {\em On the inverse spectral problem for the quasi-periodic Schr\"odinger equation.} Publ.\ Math.\ Inst.\ Hautes \'Etudes Sci.\ \textbf{119} (2014), 217--401.

\bibitem[DGL1]{DGL1} Damanik, D., Goldstein, M., Lukic, M. {\em A multi-scale analysis scheme on Abelian groups with an application to operators dual to Hill's equation.} In preparation.

\bibitem[DGL2]{DGL2} Damanik, D., Goldstein, M., Lukic, M. {\em The isospectral torus of quasi-periodic Schr\"odinger operators via periodic approximations.} In preparation.

\bibitem[DY]{DY} Damanik, D., Yuditskii, P. {\em Counterexamples to the Kotani-Last conjecture for continuum Schr\"odinger operators via character-automorphic Hardy spaces.} Preprint (arXiv:1405.6343).

\bibitem[E]{E} Eliasson, H. {\em Floquet solutions for the $1$--dimensional quasiperiodic Schr\"odinger equation.} Commun.\ Math.\ Phys.~{\bf 146} (1992), 447--482.

\bibitem[GY]{GY} Gesztesy, F., Yuditskii, P. {\em Spectral properties of a class of reflectionless Schr\"odinger operators.} J.\ Func.\ Anal.\ \textbf{351} (1999), 619--646.

\bibitem[HS]{HS} Helffer, B., Sj\"ostrand, J. {\em Analyse semi-classique pour l'\'equation de Harper} ({\em avec application \'a l'\'equation de Schr\"odinger avec champ magn\'etique}). M\'em.\ Soc.\ Math.\ France \textbf{34} (1988), 113 pp.\ (1989).

\bibitem[JM]{JM} Jitomirskaya, S., Marx, C. {\em Dynamics and spectral theory of quasi-periodic Schr\"odinger type operators.} In preparation.

\bibitem[K1]{K84} Kotani, S. {\em Ljapunov indices determine absolutely continuous spectra of stationary random one-dimensional Schr\"odinger operators.} Stochastic analysis (Katata/Kyoto, 1982), 225--247, North-Holland Math.\ Library \textbf{32}, North-Holland, Amsterdam, 1984.

\bibitem[K2]{K97} Kotani, S. {\em Generalized Floquet theory for stationary Schr\"odinger operators in one dimension.} Chaos Solitons Fractals \textbf{8} (1997), 1817--1854.

\bibitem[PR]{PR} Poltoratski, A., Remling, C. {\em Reflectionless Herglotz functions and Jacobi matrices.} Commun.\ Math.\ Phys.\ \textbf{288} (2009), 1007--1021.

\bibitem[R1]{R07} Remling, C. {\em The absolutely continuous spectrum of one-dimensional Schr\"odinger operators.} Math.\ Phys.\ Anal.\ Geom.\ \textbf{10} (2007), 359--373.

\bibitem[R2]{R11} Remling, C. {\em The absolutely continuous spectrum of Jacobi matrices.} Ann.\ of Math.\ \textbf{174} (2011), 125--171.

\bibitem[S]{Simon} Simon, B. {\em Szeg\H{o}'s Theorem and its Descendants. Spectral Theory for $L^2$ Perturbations of Orthogonal Polynomials.} M.~B.~Porter Lectures. Princeton University Press, Princeton, NJ, 2011.

\bibitem[SY1]{SY} Sodin, M., Yuditskii, P. {\em Almost periodic Sturm-Liouville operators with Cantor homogeneous spectrum.} Comment.\ Math.\ Helv.\ \textbf{70} (1995), 639--658.

\bibitem[SY2]{SY2} Sodin, M., Yuditskii, P. {\em Almost periodic Jacobi matrices with homogeneous spectrum, infinite-dimensional Jacobi inversion, and Hardy spaces of character-automorphic functions.} J.\ Geom.\ Anal.\ \textbf{7} (1997), 387--435.

\bibitem[VY]{VY} Volberg, A., Yuditskii, P. {\em Kotani-Last problem and Hardy spaces on surfaces of Widom type.} To appear in Invent.\ Math.\ (arXiv:1210.7069).

\end{thebibliography}
\end{document}